\documentclass[11pt,a4paper]{article} 
\usepackage{fullpage}
\usepackage{epsf}
\usepackage[dvips]{epsfig}
\usepackage{amssymb}
\usepackage{stmaryrd}
\usepackage[T1]{fontenc}

\usepackage{amssymb}
\usepackage{amsbsy}
\usepackage{amsfonts}

\makeindex

\newtheorem{thm}{Theorem}[section]

\newtheorem{def-not}{Definition and Notation}[section]
\newtheorem{prop}{Proposition}[section]
\newtheorem{lem}{Lemma}[section]

\newtheorem{cor}{Corollary}[section]

\newcommand{\eps}{\varepsilon}
\newcommand{\R}{\hbox{I\hskip -2pt R}}

\newcommand{\proof}{{\bf Proof :}}

\usepackage{amssymb,amsmath,amscd}

\def\proof#1{\noindent\trivlist
\item[\hskip\labelsep{\bf #1.}]
\ignorespaces}
\def\endproof{\hfill $\Box\qquad$\endtrivlist}

\def\ee{\mathrm{e}}

\begin{document}
\title{ { \bf{ Asymptotic hitting time for a simple evolutionary  model of protein folding} }}

\author{Véronique Ladret
\\
\\Laboratoire LaPCS, U.F.R. de Math\'ematiques,
\\Universit\'e Claude Bernard - Lyon 1, 50, avenue Tony Garnier,
\\B\^atiment RECHERCHE [B], Domaine de Gerland
\\69367 Lyon Cedex 07, France
\\email : veronique.ladret@univ-lyon1.fr}

\maketitle

\begin{abstract}
We consider two versions of a simple evolutionary algorithm model for protein folding at temperature zero:  the $(1+1)$-EA on the LeadingOnes problem. In this schematic model, the structure of the protein, which is encoded as a bit-string of length n, is evolved to its native conformation through a stochastic pathway of sequential contact bindings.

We study the asymptotic behavior of the hitting time, in the mean case scenario, under two different mutations: the one flip which flips  a unique bit chosen uniformly at random in the bit-string, and the Bernoulli flip which flips each bit in the bit-string independently with probability $c/n$.
 For each algorithm we prove a law of large numbers, a central limit theorem and compare the performance of the two models.
\end{abstract}
{\footnotesize{ A.M.S. S{\scriptsize{UBJECT} } C{\scriptsize{LASSIFICATION} } }}: 60J10, 60F05,  92D20, 92C05.

\section{Introduction}
Evolutionary algorithms (EAs) are  adaptive heuristic search algorithms. They  are based on the mechanisms  of natural selection and are widely used in a great variety of problems, for instance population genetics, machine learning and optimization. 
The task of the EA is to search a fitness landscape for maximal values.
A population of individuals, considered as candidate solutions to the given problem, is evolved  under steps of  mutation and steps of selection. Each individual receives a numerical evaluation, called its fitness score.  The dynamics of the EA simulates, supposedly like in natural systems, the survival of the fittest among the individuals. Thus, individuals of maximum fitness are sought.

Despite their numerous heuristic successes, mathematical results describing the behaviour of EAs are rather sparse. Among the exceptions are R. Cerf \cite{Cerf1}, \cite{Cerf2},
Y. Rabinovich and A. Wigderson \cite{RaWi},
G. Rudolph \cite{Rudolph} ,
C. Mazza and D. Piau
 \cite{MazPiau}, P. Del Moral and A. Guionnet \cite{DelGui}, J. B\'erard \cite{Berard} and J. B\'erard and A. Bienven\"ue \cite{BerBie1}, \cite{BerBie2}.

 Since EAs usually exhibit complicated dynamics, complexity results are difficult to reach and it is a common approach to consider simplified cases. Among the simplified EAs are the so called $(1+1)$-EAs.
 These are studied by H. Muhlenbein \cite{Muh}, T. B\"ack \cite{Back}, G. Rudolph \cite{Rudolph}, J. Garnier et al. \cite{GarKall}, and S. Droste et al. \cite{DJW0}, \cite{DJW1}. 
In this paper,  we study the time of convergence of two versions of a specific $(1+1)$-EA, namely the $(1+1)$-EA on the LeadingOnes problem. One of the main  motivation for studying these algorithms is that they can be used as a simple models for the protein-folding problem. Indeed, the $(1+1)$-EAs we  focus on, directly fit to the model of protein-structure prediction at temperature zero, proposed by  biophysicists A. Bakk et al. \cite{Hansen}.


\subsection{The physical model}

Proteins typically fold to a unique native or biologically active conformation on time scales from $10^{-3}$s. to 1s. However, if the dynamics of the folding process would follow a random search in the conformation space it would result in astronomical time scales. This paradox is known under the name of Levinthal's paradox \cite{Levinthal}. How do proteins fold to their native state? This is one of the intriguing problems of biophysics. Anfinsen \cite{Anfinsen} showed that the native state is genetically as well as thermodynamically determined, i.e. it corresponds to the conformation in which Gibbs free energy of the whole system is lowest. \\
There are many  hypotheses concerning the transition state (TS). One of the view is that the (TS)-dynamics consists in a pathway which carries the polypeptide (protein) to the native state through a guided descent along the Gibbs free energy landscape (J. A. Schellman \cite{Schellman},K. A. Dill et al. \cite{Dill}).\\
The protein-like model proposed by A. Bakk et al. \cite{Hansen} can be described as follows. The polypeptide chain is equipped with $n$ contact points $c_1,...,c_n$ that we will also call nodes. For $i$ from $1$ to $n$,  $c_i$ is assigned a binary contact variable $\phi_{i}$ that indicates whether it is folded ($\phi_{i}=1$) or unfolded ($\phi_{i}=0$). 
In consequence the conformation of the protein is entirely determined by the bit-string of length n,  $\phi=(\phi_1, \phi_2,...,\phi_n)$ and the native state corresponds to the bit-string where  $\phi_{i}=1, \forall 1 \leq i \leq n$, that is $(1,...,1)$. There is a bijective mapping from the conformation space onto $\{0,1\}^n$.\\

Let $i_0$ denote the smalest $i \in \{1,...,n\}$ for which $c_i$ is unfolded, i.e. for which $\phi_i=0$. We call the open part of the protein the set of contact points $\{c_{i_0}, c_{i_{0}+1},...,c_n\}$. 

The assumption about the dynamics of the folding process is that each individual node is assigned an energy of $-\eps_{0}$ if $i<i_0$, zero otherwise. It can be  implemented trhough the Hamiltonian
$$ 
{\cal{H}}= -\eps_{0} ( \phi_1+ \phi_1\phi_2+...+\phi_1...\phi_n )
$$
This Hamiltonian can be re-writen in terms of the  LeadingOnes function $L$, defined on the space of conformations $\{0,1\}^n$, which counts the length of the longest prefix of ones in the bit-string:
$$
L(x)=\max\{k \geq 1:\forall 1 \leq i \leq k, x_i=1\}\cup\{0\}.
$$
Indeed,
\begin{equation}\label{fitness}
{\cal{H}}= -\eps_{0} L(\phi)=-\eps_{0}(i_0-1)
\end{equation}

In this model, there is no energy associated to the open part of the protein. It is also known under the name of the ``zipper-model''. In fact, a descent along the energy landscape means both the folding of the leftmost uncorrectly folded substructure, i.e node $c_{i_0}$, and the ``status-quo'' for the correctly folded substructures preceding it (on the left), i.e. nodes $c_i$ for $i<i_0$.  It means that the folding events occur in a specific order: they behave like the individual locks in a zipper.  
 In  the bit-string framework, lowering the Gibbs free-energy is exactly increasing the size of the longest prefix of ones.

The algorithm proposed by A. Bakk et al. \cite{Hansen} to search the state space $\{0,1\}^n$ for conformations of lowest energy, i.e. the native state,
  is based on the Monte Carlo Metropolis (MCM) method (~Binder \cite{Binder} ).  Let $T$ denote the temperature of the whole system, $k$ the usual  Boltzmann constant, and put $\beta =1/kT$. The algorithm proceeds iteratively as follows. The individual (bit-string) at time $k$, $X_k$, undergoes a mutation to a new conformation $X_{k}'$, through a stochastic process that will be described later. Now, $X_{k}'$ is selected to form the new individual at time $k+1$, $X_{k+1}$, with probability
\begin{eqnarray*}
 P_{\mbox{accept}} & = &\min\left( 1,\exp(-\beta\Delta{\cal{H}} ) \right), \mbox{ with } \\
\Delta{\cal{H}}  & = &{\cal{H}}(X_{k}')-{\cal{H}}(X_k).
\end{eqnarray*}
Otherwise $X_{k+1}$ is a repeat of the old configuration $X_{k}$.\\

In this paper, we  concentrate on  the MCM  model taken at temperature zero, that we  denote by  $ \mbox{MCM}_{0} $, as well as on a very close version of this algorithm. These algorithms are directly connected to the $(1+1)$-EA on the LeadingOnes problem. We recall that the dynamics of $(1+1)$-EAs can be formalized through discrete Markov chains as follows:\\ 

\subsection{The $(1+1)$-EA approach to native conformation prediction}

In the protein model we wish to minimize the Hamiltonian ${\cal{H}} $, which is equivalent, according to (\ref{fitness}), to maximize the LeadingOnes function. More generally, the  goal of $(1+1)$-EAs is to optimise some fitness function $f: {\{0,1\}}^n \rightarrow {\R}$. The algorithm proceeds as follows: a  unique individual, or bit-string,  is evolved under the following two-steps iterative process:\\
1.{\bf{Mutation:}} \\
As in the MCM method, at every evolutionary step, known as generation, the individual in the current population at time $k$, $X_k$, undergoes a random walk to a new individual $X_{k}'$.\\
{\bf{ 2. Selection:}}\\
$X_k$ and $X_{k}'$ are evaluated in terms of their fitness value. Then,
the one with the highest fitness score  is selected to form the generation at time $k+1$, $X_{k+1}$:
\begin{equation}\label{rule1}
\mbox{If } f(X_{k}')> f(X_k),\mbox{ then } X_{k+1}=X'_k. \mbox{ Otherwise }  X_{k+1}=X_k.
\end{equation}

The notation $(1+1)$-EA accounts for the fact that we select among one parent and one child. 
   Here, we  focus on the mean case scenario, in which the first individual $X_0$ is chosen uniformly at random in $\{0,1\}^n$. The  reason for this restriction  is that the mean case is  easier to manage, on a mathematical level, in the LeadingOnes framework. When  the fitness function is precisely $L$, this algorithm will be denoted by 
 $(1+1)_{ \mbox{\tiny{$ L$}} } $-EA. We  notice that in the case of a fitness landscape with local maxima the $(1+1)$-EA method could end in a suboptimal search. But, as long as we consider the LeadingOnes problem whose fitness has no local maxima, we are not worried about that.
In the  litterature, there is no actual consensus, in the definition of $(1+1)$-EAs, on the selection rule. It is sometimes taken to be the following sligthly different one:

\begin{equation}\label{rule2}
\mbox{If } f(X_{k}')\geq f(X_k),\mbox{ then } X_{k+1}=X'_k. \mbox{ Otherwise }  X_{k+1}=X_k.
\end{equation}
For example, Garnier et al. \cite{GarKall} consider the first version of  the selection rule (\ref{rule1}),  whereas Droste et al. \cite{DJW0}, \cite{DJW1}  focus the second version  (\ref{rule2}).

 We notice that  in the LeadingOnes framework, the variant of $(1+1)_{ \mbox{\tiny{$ L$}} }$-EA with  mutation rule (\ref{rule2}) is nothing else but  $\mbox{MCM}_0$. In order to discriminate them, we shall keep the name  $\mbox{MCM}_0$.  The only difference between these two very close algorithms being  that $ \mbox{MCM}_{0} $ accepts candidates $X_k'$ whose energy is the same as the one of $X_k$, whereas $(1+1)_{ \mbox{\tiny{$ L$}} } $-EA does not.


\subsection{Statement of the results} 

 $T_n$ (respectively $\widehat{T}_n$) denotes the  hitting time until some optimal (with regard to the fitness function) conformation or individual is sampled by the $(1+1)$-EA (by $\mbox{MCM}_0$ respectively).  We focus on both $\mbox{MCM}_0$ and $(1+1)_{ \mbox{\tiny{$ L$}} } $-EA in the mean case scenario, 
under two different kinds of mutation: 
the  one flip,  which flips a unique bit chosen uniformly at random in the bit-string, and the Bernoulli flip, which flips each bit in the bit-string independently with probability $c/n$.  
As briefly recalled in the introduction, there has already been some work about the complexity of some  $(1+1)$-EAs:\\
 From Droste, Jansen and Wegener \cite{DJW0}, $E(\widehat{T}_n)=\Theta(n\ln n)$ for the Bernoulli flip applied to a linear fitness function on $\{0,1\}^n$. 
From Droste et al.\ \cite{DJW1}, the LeadingOnes function, is solvable in mean time $\Theta(n^2)$ in the Bernoulli flip scenario.

Garnier et al. \cite{GarKall} study the OneMax function $|\cdot|$, which counts the number of ones in the bit-string, in the one flip and the Bernoulli flip frameworks, that is
$$
|x|=\sum_{i=1}^nx_i.
$$
In the one flip case,  $(T_n-n\ln n)/n$ converges in distribution to $ -\ln 2-\ln Z$. In the  Bernoulli flip case, 
$(T_n-c^{-1}e^cn\ln n)/n$ converges in distribution to $ -c^{-1}e^c\ln Z+C(c)$, where the law of $Z$ is exponential of parameter $1$  and $C(c)$ is some $c$-dependent constant.
\\
We prove the analog for the LeadingOnes problem of the result of Garnier et al.\ \cite{GarKall}. This improves on the result of Droste et al.\ \cite{DJW1}. We prove a law of large numbers, a central limit theorem, and we compare the performance of the two models. Finally, we prove that the distribution of the hitting time of $ \mbox{MCM}_{0} $,  $\widehat{T}_n$, is the same as the one of $T_n$ in both the one flip and the Bernoulli flip scenari.

\begin{thm}[one flip case] 
\label{t.1flip}
(i) For $n \geq 1$, 
$
E(T_n)=n^2/2$;
\\
(ii)  As  $n\to\infty$,
$T_n/E(T_n)$ converges in probability to $1$.
\\
(iii)
 As  $n\to\infty$,
$(T_n-E(T_n))/n^{3/2}$
converges in distribution to a centered 
Gaussian random variable of variance $3/4$.
\end{thm}
\begin{thm}[Bernoulli flip case] 
\label{t.bflip}
i)
As  $n\to\infty$,
$E(T_n)\sim m(c)\,n^2$, with
$$
m(c):=(e^c-1)/(2c^2).
$$
ii) As  $n\to\infty$,
$T_n/E(T_n)$ 
converges in probability to $1$.
\\
iii)
Furthermore,
$(T_n-m(c)\,n^2)/n^{3/2}$
converges in distribution to a centered 
Gaussian random variable of variance $\sigma^2(c)$,
with
$$
\sigma^2(c):=3(e^{2c}-1)/(8c^3).
$$
\end{thm}
Note that $m(c)>1/2$ for every $c>0$. 
\begin{cor}
\label{c.comp}
As $n\to\infty$, 
$E(T_n)$ for the Bernoulli flip case is greater than
$E(T_n)$ for the one flip case, for any value of $c$.
\end{cor}
\begin{thm} \label{t3_MCM}
In both the one flip  and the bernoulli flip case,
$T_n$ and $\widehat{T}_n$ have the same distribution.
\end{thm}

\section{Proof of Theorem~\ref{t.1flip}}

The law of $T_n$, conditioned by $|X_0|$, is the law of 
a sum of geometric random variables, see Lemma~\ref{l2.1}.
This yields Part~(i) of the theorem.
Since the CLT implies the law of large numbers, we then prove 
the CLT of Part~(iii). 

The distribution of $T_n$ stems from a simple observation.
In the $(1+1)$-EA on the LeadingOnes problem, 
a mutation is accepted if and only if it adds $1$ to 
the number of leading ones.
As a consequence, in the one flip framework,  
the chain jumps when the leftmost zero is flipped.
The other flips leave the chain unchanged.
Thus, the zeroes in $X_0$ are successively flipped, 
from left to right, until one hits the optimal individal (conformation) $(1,1,\ldots,1)$.

Here and later in the paper, $\eps=1/n$ when we deal 
with algorithms on strings of length $n$, 
$|x|$ is the number of ones in $x$,
the geometric law ${\cal G}(p)$ of parameter $p$ is defined by
$$
{\cal G}(p)=\sum_{n\ge1}p\,(1-p)^{n-1}\delta_n,
$$
and the negative binomial law
of parameter $(k_0,p)$, ${\cal{NB}}(k_0,p)$, puts the following mass on $k\ge k_0$:
$$
{k-1\choose k-k_0}p^{k_0}(1-p)^{k-k_0}.
$$
\begin{lem}\label{l2.1}
If $|X_0|=n-k_0$, 
$T_n$ is the sum of $k_0$ i.i.d.\  ${\cal G}(\eps)$ random variables.
Thus, the law of  $T_n$ is negative binomial of parameter $(k_0,\eps)$.
\end{lem}
\proof{Proof of Lemma~\ref{l2.1}}
Let $\tau_0=0$ and, for every $k\ge0$,
\begin{equation}\label{e4}
\tau_{k+1}=\inf\{i \geq \tau_{k}\,; X_i \neq X_{\tau_{k}}\},
\quad
\sigma_{k+1}= \tau_{k+1}-\tau_{k},
\quad
\widetilde{X}_k=X_{\tau_k}.
\end{equation}
In words, ${\widetilde X_k}$ denotes the position of the chain after 
its $k$-th jump, and $|{\widetilde X_k}|=n-k_0+k$. 
The leftmost zero of ${\widetilde X_k}$ is flipped after a sojourn at ${\widetilde X_k}$ of length $\sigma_{k+1}$.
Thus, $(\sigma_k)_k$ is i.i.d.\ of law ${\cal G}(\eps)$. 
It remains to note that
$$
T_n=\sigma_1+\cdots+\sigma_{k_0}.
$$
\endproof
\proof{Proof of Part~(iii)}
Let $E_k$ denote the conditioning on $\{|X_0|=n-k\}$.
Since $X_0$ is uniform, the law $\mu$ of $|X_0|$ is binomial $(n,1/2)$.
From Lemma~\ref{l2.1},
under $P_k$, $T_n$  is the sum of $k$ i.i.d.\ geometric random variables of parameter $\eps$.
Thus,
\begin{equation}
\label{e.ek}
E_{k}(\ee^{-\alpha T_n})=\ee^{-\alpha k}\,\eps^{k}\left[1-(1-\eps)\,\ee^{-\alpha}\right]^{-k}.
\end{equation}
Set $\Theta_n=(T_n-n^2/2)/n^{3/2}$.
The decomposition of the Laplace transform of $\Theta_n$
along the values of $|X_0|$, the explicit form of $\mu$, and Equation~(\ref{e.ek}) yield together that
$$
E(\ee^{-\alpha \Theta_n})
=
\ee^{\sqrt{n}\alpha/2}\,2^{-n}\sum_{k=0}^n{n\choose k}\,\ee^{-\alpha k/n^{3/2}}\,\eps^{k}\left[1-(1-\eps)\,\ee^{-\alpha/n^{3/2}}\right]^{-k}.
$$
This can be rewritten as
$$
E(\ee^{-\alpha \Theta_n})
=
\ee^{\sqrt{n}\alpha/2}\,2^{-n}\,(1+\beta_n)^{n},
$$
with
$$
\beta_n=\eps\,\ee^{-\alpha/n^{3/2}}\left[1-(1-\eps)\,\ee^{-\alpha/n^{3/2}}\right]^{-1}.
$$
Recall that $\eps=1/n$. The expansion of $\beta_n$ reads
$$
\beta_n=1-\alpha/\sqrt{n}+\alpha^2/n+o(1/n).
$$
This implies that
$$
E(\ee^{-\alpha \Theta_n})
\to
\ee^{3\alpha^2/8},
$$
as $n\to\infty$. This concludes the proof.
\endproof


\section{Proof of Theorem~\ref{t.bflip}}

We first describe the law of $T_n$ conditionally on the values taken 
by $L$ along the path of $(X_k)_k$ before $T_n$, that is, 
until the optimal individual $(1,1,\ldots,1)$ is hit.
This law is the law of a sum of independent geometric random variables, see Lemma \ref{Le3.3}. We deduce the overall law of $T_n$, see 
Proposition~\ref{pro3.1}. This yields part (i) of the theorem.

As in the previous section, since CLT implies the law of large numbers, we then prove the CLT of part (iii).

We recall that the $(1+1)$ EA, in the LeadingOnes framework, accepts a mutation if and only if the number of leading ones is increased. Hence the dynamics of the Bernoulli flip algorithm proceeds as follows: the chain  jumps to a new individual, at time $k+1$, if and only if the leading ones of $X_k$ are left unchanged and  its leftmost zero is flipped, no matter which values are taken by the other bits.

Here and later in the paper, $\eps=c/n$ when we deal with algorithms on strings of length $n$.
For all $i\geq 0$, let  $p(n,i)=\eps(1-\eps)^i$.
As in the one flip framework, $ \widetilde{X}_k$ denotes the position of the chain after its $k$-th jump. We also keep the same definition for $\sigma_k$ and $\tau_k$.\\
For all $k \geq 0$, let $$\ell_k=L(\widetilde{X}_k).
$$
Let $Y_0$ be such that $$X_0=(1^{\ell_0},0,Y_0).$$
For $k \geq 1$, define $Y_k$ and $W_k$ by 
$$
\widetilde{X}_k=(1^{\ell_k-1},1,W_k)=(1^{\ell_k},0,Y_k).
$$
Lemmas~\ref{lm} and \ref{Lo} below are needed to compute the law of 
of $T_n$ conditionally on $(\ell_j)$ in Lemma~\ref{Le3.3}.

\begin{lem}\label{lm}
(i) For all $k \geq 0$,
${\sigma}_k$ depends on the past only through the last score $\ell_{k-1}$.
That is, the law of $\sigma_k$, conditionally on $(X_t)_{\{t<\tau_k\}}$, is the law of $\sigma_k$, conditionally on $\ell_{k-1}$.
\\
(ii)
For all $k \geq 0$, given $\{\ell_{k-1}=i\}$, the law of $\sigma_k$  is $\mathcal{G}(p(n,i))$.
\end{lem}

\proof{Proof}
The sojourn time ${\sigma}_k$ is the time the algorithm takes to jump from $\widetilde{X}_{k-1}$ to $\widetilde{X}_{k}$.
The leftmost zero of $\widetilde{X}_{k-1}$ is in position $\ell_{k-1}+1$.
Thus, one needs to flip the  $(\ell_{k-1}+1)$-th bit, while leaving the first $\ell_{k-1}$ bits unchanged.Thus,
$$
P({\sigma}_k=t\,|\,\widetilde{X}_0,\ldots,\widetilde{X}_{k-1} )=\eps(1-\eps)^{\ell_{k-1}}\left[1-\eps(1-\eps)^{ \ell_{k-1} } \right]^{(t-1)}.
$$
\endproof
\begin{lem}\label{Lo}
For all $k \geq 1$ let ${\widehat {\cal F}}_k={\sigma}\{\sigma_i, \ell_j: i \leq k, j \leq k-1\}$, then
the law of $\ell_{k}$ conditionally on ${\widehat {\cal F}}_k$ is the law of  $\ell_{k}$  conditionally on $\ell_{k-1}$.

\end{lem}

\proof{Proof}
Let $P_{i}$ denote the probability given $\{\ell_0=i\}$.
Using the strong Markov property:
\begin{equation}\label{e6}
P\left(\ell_k=i_k | \sigma_k=t_k, \ell_{k-1}=i_{k-1},\ldots,\sigma_1=i_1, \ell_0=i_0\right)=
P\left(\ell_1=i_k | \sigma_1=t_k, \ell_{0}=i_{k-1}\right)
\end{equation}
Since $\{ \ell_1=i_1 , \sigma_1=t, \ell_0=i_0\}=\{ L({ X}_t)=i_1, L(X_{t-1})=\ldots=L(X_0)=i_0\}$, since $\{\sigma_1=t, \ell_0=i_0\}=\{ L({ X}_t) \neq i_0, L(X_{t-1})=\ldots=L(X_0)=i_0\}$ and using the Markov property on $\left(L(X_t)_t\right)$  we derive the following expression: 

$$
 P\left(\ell_1=i_k | \sigma_1=t_{k}, \ell_0=i_{k-1}\right)=
 \frac{ P\left(L( X_1)=i_k|L( X_{0})=i_{k-1}\right)}{ P\left(L( X_1) \neq i_{k-1}|L( X_{0})=i_{k-1}\right)}
$$
This quantity is indepedent from $t_k$, hence if we reconsider (\ref{e6}):
$$
P\left(\ell_k=i_k | \sigma_k=t_k, \ell_{k-1}=i_{k-1},\ldots,\sigma_1=i_1, \ell_0=i_0\right)=
P\left(\ell_k=i_k | \ell_{k-1}=i_{k-1}\right)
$$
\endproof

\begin{lem}\label{Le3.3}
Conditionally on $\{\ell_0=i_0,\ldots,\ell_{J-1}=i_{J-1},\ell_J=n\}$, $T_n$ is the sum of $J$ independent geometric random variables with respective parameters $p(n,i_0),\ldots,p(n,i_{J-1})$.
\end{lem}
\proof{Proof}
Given the successive LeadingOnes scores $\{\ell_0=i_0,\ldots,\ell_J=n\}$ until the optimal individual is hit,  $T_n=\sum_{k=1}^J\sigma_k$. Thus,
\begin{equation}\label{e7}
P(T_n=t|\ell_J=n,\ldots,\ell_0=i_0)
=\sum_{t_1+\cdots+t_k=t}P(\sigma_J=t_J,\ldots, \sigma_1=t_1|\ell_J=n,\ldots,\ell_0=i_0)\\
\end{equation}
Using Lemmas ( \ref {Lo}) and ( \ref {lm}) we can derive by induction:
$$
P\left(\ell_J=n, \sigma_{J}=t_J,\ldots,\sigma_1=t_1,\ell_0=i_0\right)
=\prod_{k=1}^JP(\ell_k=i_k|\ell_{k-1}=i_{k-1}) P(\sigma_{k}=t_k|\ell_{k-1}=i_{k-1})
$$
Hence, since $\prod_{k=1}^JP(\ell_k=i_k|\ell_{k-1}=i_{k-1})= P(\ell_J=n,\ldots,\ell_0=i_0)$,

$$
 P(\sigma_1=t_1,\ldots,\sigma_J=t_J|\ell_{J}=n,\ldots,\ell_0=i_0)
=\prod_{k=1}^J P\left( \sigma_{k}=t_k|\ell_{k-1}=i_{k-1} \right)
$$

We can now replace this last equation in ($\ref{e7}$),

$$
P(T_n=t|\ell_J=n,\ldots,\ell_0=i_0)
=\sum_{t_1+\cdots+t_k=t}\prod_{k=1}^J P\left(\sigma_{k}=t_k|\ell_{k-1}=i_{k-1}\right)
$$
Let $q(n,i_k)$ denote the probabilty distribution of ${\cal G}(p(n,i_k))$.\\
Then using Lemma {\ref {lm}}, we can write, as we recognize a product of convolution:
$$
P(T_n=t|\ell_J=n,\ldots,\ell_0=i_0) = q(n,i_{J-1})*\cdots*q(n,i_0)(t)
$$
Thus, given that the search jumps $J$ times until the target $(1,1,\ldots,1)$ is hit and given\\$\{\ell_0=i_0,\ldots,\ell_J=n\}$, $T_n$ follows the same distribution as the sum of $J$ independent ${\cal G}(p(n,i_0)),$
$\ldots, {\cal G}(p(n,i_{J-1}))$ random variables.\\
~\\
Let us focus on $P(\ell_0=i_0,\ldots,\ell_{J-1}=i_{J-1},\ell_J=n)$.
In order to compute this quantity, we need the following Lemma:\\
\begin{lem}\label{Le3.4} Let $k \geq 1$. 
If $X_0$ is chosen uniformly in ${\{0,1\}}^n$, then, 
given $\{\ell_{k-1}=i\}$, $W_k$ follows the uniform distribution on  ${\{0,1\}}^{n-i-1}$:
$$ 
{\cal L}(W_k|\ell_{k-1}=i) ={\cal U}(\{0,1\}^{n-i-1}).
$$ 
\end{lem}

\proof{Proof}
Let us focus on the case where $k=1$.
Let $P_i$ denote the probability conditionally on $\{\ell_0=i\}$.
Let $\mu$ denote the probabilty distribution of $Y_0$ given $\{\ell_0=i\}$.
As $X_0$ is chosen uniformly in ${\{0,1\} }^n$, $\mu$ is the uniform distribution on ${\{0,1\} }^{n-i-1}$.

\begin{equation}\label{e8}
P_i(W_1=w)=\sum_{t \geq 1  }P_i(\widetilde{X}_1=(1^i,1,w), \sigma_1=t)
\end{equation}
Since $\{ X_{\sigma_1}=(1^i,1,w), {\sigma_1}=t, L(X_0)=i \}=\{  X_t=(1^i,1,w), L(X_{t-1})=\ldots=L(X_0)=i  \}$ and using the Markov property, Equation (\ref{e8}) can be rewritten:
$$
P_i(W_1=w)
=\sum_{ t \geq 1 }P_i(X_1=(1^i,1,w)){ P_i( L(X_1)=L(X_0)) }^{t-1}
$$
Hence,
\begin{equation}\label{e9}
P_i(W_1=w)
=\frac{P_i(X_1=(1^i,1,w))}{P_i(L(X_1) > L(X_0) )}
\end{equation}

We recall that the Markov chain jumps from  $X_0$ to a conformation of higher fitness at time $1$ if both none of the $\ell_0$ first ones of $X_0$ are flipped and the leftmost zero of $X_0$, is flipped. Hence,
\begin{equation}\label{e10}
P_i(L(X_1) > L(X_0) )= \eps(1-\eps)^i
\end{equation}
On the other hand, as $P_i(Y_0=u)=\mu(u)=1/2^{n-i-1}$,
\begin{equation}\label{e11}
 P_i(X_1=(1^i,1,w))=1/2^{n-i-1}\sum_{u \in { \{0,1\} }^{n-i-1} }P(X_1=(1^i,1,w)|X_0=(1^ i,0,u))
\end{equation}
If $d(w,u)$ denotes the Hamming distance between $w$ and $u$,
\begin{equation}\label{e12} 
P(X_1=(1^i,1,w)| X_{0}=(1^i,0,u))=\eps(1-\eps)^i\eps^{d(w,u)}(1-\eps)^{n-i-1-d(w,u)} 
\end{equation}

Finally Equations (\ref{e9}), (\ref{e10}), (\ref{e11}) and (\ref{e12}) together with
$$
 \sum_{u \in { \{0,1\} }^{n-i-1} }\eps^{d(w,u)}(1-\eps)^{n-i-1-d(w,u)}=1 \mbox{, }$$ yield that: 

\begin{equation}\label{ww1}
P_i(W_1=w)=1/2^{n-i-1}
\end{equation}

Now, using the strong Markov property we can derive the proof for any $k \geq 2$.
\endproof

\begin{lem}\label{lmL}If $X_0$ is chosen uniformly at random in the state space,
for all $k \geq 1$ ,
the  conditional distribution of $ \ell_k$, given $\{\ell_{k-1}=i_{k-1}\}$, satisfies:
\begin{eqnarray*}
 P\left( \ell_k=j_k|  \ell_{k-1} = j_{k-1}\right)& = & 2^{-( j_k-j_{k-1} ) } \mbox{ \hspace{0.4cm}  if \hspace{0.1cm} } i_0+1 \leq j_k < n \\
& = & 2^{-( n-j_{k-1}-1 ) } \mbox{ \hspace{0.2cm}  if \hspace{0.1cm} } j_k = n\\
\end{eqnarray*}
\end{lem}
\proof{Proof}
This is a direct consequence of Lemma \ref{Le3.4}
\endproof

 Now that we know the probability distribution of  the sequence of the successive LeadingOnes scores until the target individual $(1,1,\ldots,1)$ is hit as well as the distribution of $T_n$ conditional to the values taken by these LeadingOnes scores, we can compute the probability distribution of $T_n$:

\begin{prop} \label{pro3.1}
If $X_0$ is chosen uniformly at random in ${\{0,1\}}^n$,
then, the probability distribution of $T_n$ satisfies:
\begin{equation}\label{dis}
P(T_n=t) 
= 
\frac{1}{2^n}\sum_{J}\sum_{\{o \leq i_0 \leq i_1 \leq \ldots \leq i_J=n\}}q(n,i_{J-1})*\cdots*q(n,i_0)(t)
\end{equation}
\end{prop}

\proof{Proof}
$X_0$ being chosen uniformly at random in the search space, $P(\ell_0=i_0)=\frac{1}{2^{i_0+1}}$.
Thus, applying Lemma \ref{lmL},
$$
P(\ell_0=  i_0,\ldots, \ell_J=n)= \frac{1}{ 2^n}
$$
The result is then a direct consequence of Lemma \ref{Le3.3}. 
\endproof


\proof{Proof of part (iii)}
Set ${\Theta}_n=(T_n -\frac{n^2(e^c-1)}{2c^2})|n^{3/2}$.

Thus,
\begin{equation}\label{e15}
E(\exp(-\alpha {\Theta}_n))=\exp\left(\alpha\frac{\sqrt{n}}{2c^2}(e^c-1)\right)E\left(\exp(-\alpha\frac{T_n}{n^{3/2}})\right)
\end{equation}
According to the distribution of $T_n$ given by (\ref{dis}):
$$
E\left(\exp(-\alpha\frac{T_n}{n^{3/2}})\right) = 
 \frac{1}{2^n}\sum_{J}\sum_{ i_0<\cdots<i_{J-1}}E\left(\exp(-\frac{\alpha}{n^{3/2}}({\cal G}(p(n,i_0)+\cdots+{\cal G}(p(n,i_{J-1}))\right)$$
Since the variables $({\cal G}(p(n,i_k))_{i_k}$ are independent,
\begin{eqnarray*}
E\left(\exp(-\alpha\frac{T_n}{n^{3/2}})\right) 
& = & \frac{1}{2^n}\sum_{J}\sum_{ i_0<\cdots<i_{J-1}}\prod_{k=0}^{J-1}E\left(-\frac{\alpha}{n^{3/2}}{\cal G}(p(n,i_k)\right)\\
& = & \frac{1}{2^n}\prod_{k=0}^{n-1}\left(1+\phi_{k}(\frac{\alpha}{n^{3/2}})\right)\\
\end{eqnarray*}
$\mbox{ where } \phi_{k}(\alpha) \mbox { denotes the Laplace transform of }{\cal G}(p(n,i_k))$.

$$
\phi_{k}(\alpha)= \frac{ e^{-\alpha}p(n,i_k) }{ 1-(1-p(n,i_k))e^{ -\alpha}  }.
$$
Recalling that $p(n,i_k)=\eps(1-\eps)^{i}$ and $\eps=c/n$, we derive that:
$$
\prod_{k=0}^{n-1}\left(1+\phi_{k}(\frac{\alpha}{n^{3/2}})\right) \simeq_{n \rightarrow +\infty} 2^n\exp\left(-\alpha \frac{ \sqrt{n} }{2c^2}(e^c-1) \right)\exp\left( \frac{ 3{\alpha}^2 }{8c^3}\frac{ (e^{2c}-1 )}{2} \right)
$$
Thus, replacing this in ($\ref{e15}$), as $n$ goes to $\infty$:
$$
E(\exp(-\alpha{\Theta}_n ))\simeq \exp\left( \frac{ 3{\alpha}^2 }{8c^3}\frac{ (e^{2c}-1 )}{2} \right)
$$
We  recognize the Laplace transform of a centered gaussian variable of variance $ \frac{ 3(e^{2c}-1 )}{8c^3}$. It ends the proof.


\section{Proof of Theorem~\ref{t3_MCM}}

In the $\mbox{MCM}_0$ framework, the algorithm has the possibility to visit several distinct protein conformations (individuals) whose fitness has the same value before it finally jumps to a new individual with higher fitness  score. This is not allowed  in  the $(1+1)_{ \mbox{\tiny{$ L$}} }$-EA where  the individual at time $k$, $X_k$, is not allowed to jump to an individual $X_{k+1}$ with the same fitness score and that would not be $X_k$ itself.

As in the previous sections  we denote by  $(\widetilde{X}_k)_k$  the chain defined by the protein conformations taken at the times of fitness jumps $(\tau_k)_k$.\\

Now, we briefly sketch the proof of Theorem \ref{t3_MCM}.\\ 
{\bf{Proof in the one flip case.}}\\
Here, we put $\eps=1/n$.\\ 
 In  the $( 1+1)$ framework, once $X_0$ has been sampled, the path $(\widetilde{X}_k)_k$ becomes entirely deterministic: going trough it means exactly flipping, one at a time, from left to right the zeros of $X_0$. This is not true in the case of $\mbox{MCM}_0$ and we cannot adapt directly the proof of Theorem \ref{t.1flip}.

 The idea of the proof is close to the one Theorem  \ref{t.bflip}. 
First  we  consider the law of $\widehat{T}_n$ conditionally on the values  taken by $L$ along the path $(X_k)_k$ until the ground state is sampled. The same basics arguments apply. Lemma \ref{Lo} remains true and we find that, as in the one flip $(1+1)$ framework, $(\sigma)_k$ is i.i.d. of law ${\cal{G}}(\eps)$. Now, Lemma \ref{Le3.3}  can easily be adapted:

\begin{lem}\label{Le4.1}
Conditionally on $\{\ell_0=i_0,\ldots,\ell_{\mbox{\footnotesize{$J$}}-1}=i_{\mbox{\footnotesize{$J$}}-1},\ell_{\mbox{\footnotesize{$J$}}}=n\}$, $\widehat{T}_n$ is the sum of $\mbox{\footnotesize{$J$}}$ i.i.d. geometric random variables with  parameter $\eps$, i.e. $\widehat{T}_n$ is negative binomial of parameter $(\mbox{\footnotesize{$J$}}, \eps)$.
\end{lem}

Lemma \ref{Le3.4} still holds:

\proof{Proof}
The proof is a copy the one of Lemma \ref{Le3.4}  up to Equation (\ref{e9}):

$$
P_i(W_1=w)
=\frac{P_i(X_1=(1^i,1,w))}{P_i(L(X_1) > L(X_0) )}
$$
In the one flip scenario, a unique bit at a time is flipped during the step of mutation. In consequence, in order to sample $(1^i,1,w)$ at time $1$, we need to have sampled $(1^i,0,w)$ at time $0$. Thus,
\begin{equation}\label{e16}
P_i(X_1=(1^i,1,w)) = P(X_1=(1^i,1,w)|X_0=(1^ i,0,w))P_i(Y_0=w)
                   = 1/(n2^{n-i-1})
\end{equation}
On an another hand,
\begin{equation}\label{e17}
P_i(L(X_1) > L(X_0) )=1/n
\end{equation}
Equalities (\ref{e16}) and (\ref{e17}) applied to Equation  (\ref{e9})  return the result for $k=1$. Finally, the strong Markov property  ends the proof for  $k>1$.
 
Now, we can prove Theorem \ref{t3_MCM} in the one flip case:
\proof{Proof of  Theorem \ref{t3_MCM} }
From the above, we derive the distribution of $\widehat{T}_n$:
\begin{eqnarray*}\label{dis2}
P(\widehat{T}_n=t) & = & 
\frac{1}{2^n}\sum_{J}\sum_{\{o \leq i_0 \leq i_1 \leq \ldots \leq i_J=n\}}{\cal{NB}}(\mbox{\footnotesize{$J$}},\eps)(t) \\
 & = & \sum_{J} \frac{1}{2^n}{n \choose \mbox{ \footnotesize{$J$} } }{\cal{NB}}(\mbox{\footnotesize{$J$}},\eps)(t)
\end{eqnarray*}
 We recall from Lemma \ref{l2.1} that conditionally on $|X_0|=n-\mbox{\footnotesize{$J$}}$, $T_n$ is negative binomial of parameter $(\mbox{\footnotesize{$J$}},\eps)$. As $X_0$ is chosen uniformly at random in $\{0,1\}^n$, 
it yields that for all $ t \geq 0$:
 $$
P({T}_n=t)=P(\widehat{T}_n=t) 
$$
\endproof
{\bf{Proof in the Bernoulli flip case.}}\\
The dynamics of $\mbox{MCM}_0$ and  $(1+1)_{ \mbox{\tiny{$ L$}} } $ slighlty differ. Though,
we notice that in the $\mbox{MCM}_0$ framework,  Lemmas \ref {lm}, \ref{Lo},  \ref{Le3.4}, \ref{lmL} still hold. Now, since the proof of Proposition \ref{pro3.1}, is entirely based on these lemmas, we derive that the probabilty distribution of the hitting time is the same in both the $\mbox{MCM}_0$ and the $(1+1)_{ \mbox{\tiny{$ L$}} } $-EA scenari.\\
$\mbox{ }$ 
\endproof



\section{Conclusion}

After examining the two versions (one flip and Bernoulli flip) of the EAs on which we focused, we reach the following conclusion: As $ (e^c-1)/c^2 >1$ for all $c \in \R_{+}$, the expected value of the hitting time is higher in the Bernoulli flip  than in the one flip. Thus, we can conclude that the one flip  performs better than any Bernoulli flip, in terms of the expected hitting time. The same conclusion has already been derived by Garnier et al.\ \cite{GarKall} for the OneMax problem.

This better performance of the one flip suggests that, despite the ability of the Bernoulli flip to jump from any region of the search space to any other, in  a single iteration of the search process, the Bernoulli flip results in a slower convergence to a given individual, in the LeadingOnes framework.

 In order to explain this phenomenon, as the Markov chain which models our $(1+1)$ search process accepts a mutation only in case of an increase in the number of leading ones, we should point out the following facts: in the Bernoulli flip framework,  the closer the algorithm gets to the target  individual, the longer the algorithm waits until it jumps; on the other hand, in the  one flip case, the number of leading ones currently present in the bit-string  does not interfere in the distribution of the time taken for the search to jump. Also, this probability distribution remains stable as the search draws near to the optimal individual.


\begin{thebibliography}{99}






\bibitem{Anfinsen}
C.~B. Anfinsen.
\newblock Principles that govern the folding of a protein chain ({N}obel
  lecture).
\newblock {\em Science}, 191:223--230, 1973.

\bibitem{Back}
T.~B\"ack.
\newblock Optimal {M}utation {R}ate in {G}enetic {S}earch.
\newblock In S.{F}orrest, editor, {\em Proceedings of the $5^{th}$
  {I}nternational {C}onference on {G}enetic {A}lgorithms}, pages 2--8. Morgan
  Kaufmann, 1993.

\bibitem{Hansen}
A.~Bakk, J.~S. H{\o}ye, A.~Hansen, S.~Sneppen, and M.~H. Jensen.
\newblock Pathways in two-state protein folding.
\newblock {\em Biophys.J.}, 79:2722--2727, 2000.

\bibitem{Berard}
J.~B\'erard.
\newblock Genetic algorithms in random environments: two examples.
\newblock Preprint.

\bibitem{BerBie1}
J.~B{\'e}rard and A.~Bienven{\"u}e.
\newblock Convergence of a genetic algorithm with finite population.
\newblock In {\em Mathematics and computer science (Versailles, 2000)}, pages
  155--163. Birkh\"auser, Basel, 2000.

\bibitem{BerBie2}
J.~B{\'e}rard and A.~Bienven{\"u}e.
\newblock Un principe d'invariance pour un algorithme g\'en\'etique en
  population finie.
\newblock {\em C. R. Acad. Sci. Paris S\'er. I Math.}, 331(6):469--474, 2000.

\bibitem{Binder}
K.~Binder.
\newblock {\em Applications of the Monte Carlo Method in Statistical Physics}.
\newblock Springer Verlag, Berlin, 1987.

\bibitem{Cerf2}
R.~Cerf.
\newblock The dynamics of mutation-selection algorithms with large population
  sizes.
\newblock {\em Ann. Inst. H. Poincar\'e Probab. Statist.}, 32(4):455--508,
  1996.

\bibitem{Cerf1}
R.~Cerf.
\newblock Asymptotic convergence of genetic algorithms.
\newblock {\em Adv. in Appl. Probab.}, 30(2):521--550, 1998.

\bibitem{Dill}
K.A. Dill, K.~M. Feibig, and H.~S. Chan.
\newblock Cooperativity in protein folding kinetics.
\newblock {\em Proc.Natl.Acad.Sci.USA}, 90:1942--1946, 1993.

\bibitem{DJW}
S.~Droste and et~al.
\newblock On the analysis of the (1+1) evolutionary algorithm.

\bibitem{DJW1}
S.~Droste, T.~Jansen, and I.~Wegener.
\newblock On the analysis of the (1+1) evolutionary algorithm.
\newblock Technical Report CI-21/98,Univ. Dortmund, Collaborative Research
  Center 531, 1998.

\bibitem{DJW0}
S.~{D}roste, {T}. {J}ansen, and {I}. {W}egener.
\newblock A rigorous complexity analysis of the $(1+1)$ evolutionary algorithm
  for linear functions with boolean inputs.
\newblock In {\em Proceedings of the {F}ith {I}{E}{E}{E} {I}nternational
  {C}onference on {E}volutionary {C}omputation}. {I}{E}{E}{E} {P}ress, 1998.

\bibitem{GarKall}
J.~Garnier, L.~Kallel, and M.~Schoenauer.
\newblock Rigorous hitting times for binary mutations.
\newblock {\em Evolutionary Computation}, 7(2):173--203, 1999.

\bibitem{Levinthal}
C.~Levinthal.
\newblock Principles that govern the folding of a protein chain ({N}obel
  lecture).
\newblock {\em J. Chem. Phys.}, 65:44--45, 1968.

\bibitem{MazPiau}
Christian Mazza and Didier Piau.
\newblock On the effect of selection in genetic algorithms.
\newblock {\em Random Structures Algorithms}, 18(2):185--200, 2001.

\bibitem{DelGui}
P.~Del Moral and A.~Guionnet.
\newblock On the stability of interacting processes with applications to
  filtering and genetic algorithms.
\newblock {\em Ann. Inst. H. Poincar\'e Probab. Statist.}, 37(2):155--194,
  2001.

\bibitem{Muh}
H.~M\"uhlenbein.
\newblock How genetic algorithms really work: {I}. mutation and hill-climbing.
\newblock In R.{M}anner and {B}.{M}anderick, editors, {\em Proceedings of the
  $2^{nd}$ {C}onference on {P}arallel {P}roblems {S}olving from {N}ature},
  pages 15--25. Morgan Kaufmann, 1992.

\bibitem{RaWi}
Y.~Rabinovich and A.~Wigderson.
\newblock Techniques for bounding the convergence rate of genetic algorithms.
\newblock {\em Random Structures Algorithms}, 14(2):111--138, 1999.

\bibitem{Rudolph}
G.~Rudolph.
\newblock {\em Convergence {P}roperties of {E}volutionary {A}lgorithms}.
\newblock Kovac, Hamburg, 1997.

\bibitem{Schellman}
J.~A. Schellman.
\newblock The factors affecting the stability of hydrogen polypeptide
  structures in solution.
\newblock {\em J. Phys. Chem.}, 62, 1485-1494.

\end{thebibliography}
\end{document}